\documentclass[12pt]{article}
\usepackage{amsfonts}
\usepackage{hyperref}
\usepackage{mathrsfs}
\usepackage{indentfirst}
\usepackage{amsmath}
\usepackage{amssymb}
\usepackage[amsmath,thmmarks]{ntheorem}
\usepackage{cite}
\usepackage{graphicx}
\usepackage{tikz}

\newtheorem{definition}{\bf Definition}[section]
\newtheorem{lemma}{\bf Lemma}[section]
\newtheorem{theorem}{\bf Theorem}[section]

\newtheorem{proposition}{\bf Proposition}[section]

\begin{document}
\setcounter{page}{1}

\title{{\textbf{The characterization of general pseudo-homogeneous t-norms}}\thanks {Supported by
the National Natural Science Foundation of China (No. 12471440)}}
\author{Feng-qing Zhu$^{1}$\footnote{\emph{E-mail address}: zfqzfq2020@126.com}, Xue-ping Wang$^2$\footnote{Corresponding author. xpwang1@hotmail.com; fax: +86-28-84761502},\\
\emph{\small {1. School of Sciences, Southwest Petroleum University, Chengdu 610500,}}\\\emph{\small {Sichuan, People's Republic of China}}\\
\emph{\small {2. School of Mathematical Sciences, Sichuan Normal University, Chengdu 610066,}}\\
\emph{\small {Sichuan, People's Republic of China}}}

\newcommand{\pp}[2]{\frac{\partial #1}{\partial #2}}
\date{}
\maketitle
\begin{quote}
{\bf Abstract} This article first introduces the concept of a general pseudo-homogeneous triangular norm. It then gives some properties of general pseudo-homogeneous triangular norms. Finally, it characterizes all general pseudo-homogeneous triangular norms completely.

\end{quote}

\section{Introduction}
 The homogeneous functions play a central role in many fields such as image processing, decision making, etc., see, e.g., \cite{Cabrera2014,Fuss1978,Jurio2013, Jurio2014}, which directly leads the homogeneity of peculiar aggregation functions to a new hot direction with rich achievements. For example, R\"{u}ckschlossov\'{a} established a complete characterization of homogeneous aggregation functions \cite{T2005}. Alsina, Frank and Schweizer showed that a triangular norm (t-norm for short) $T$ is homogeneous of order $k>0$ if and only if $k=2$ and $T=T_P$, or $k=1$ and $T=T_M$ and that $S$ is a homogeneous t-conorm of order $k>0$ if and only if $k=1$ and $S=S_M$ \cite{Alsina2006}. In particular, Ebanks introduced a more relaxed homogeneity, called a quasi-homogeneity, which is defined as $T(\lambda x, \lambda y) =\varphi^{-1}(\psi(\lambda)\varphi(T(x, y)))$ for any $x, y, \lambda \in[0, 1]$ where $\psi :[0, 1] \rightarrow[0, 1]$ is an arbitrary function, $\varphi:[0, 1] \rightarrow[0, \infty)$  is a continuous injection and $T$ is a t-norm \cite{Ebanks1998}. In addition, Xie, Su and Liu \cite{Xie2016} gave a concept of a pseudo-homogeneous t-norm, i.e., $T$ is called a pseudo-homogeneous t-norm if, for all $x, y,\lambda\in[0, 1]$, it holds that $T(\lambda x,\lambda y)=F(\lambda,T(x,y))$ where $F:[0,1]^2\rightarrow[0,1]$ is a continuous and increasing function that fulfills the boundary condition $F(x,1)=0\Leftrightarrow x=0$, and they described the pseudo-homogeneous t-norm completely. Clearly, the pseudo-homogeneous t-norm must be continuous. Therefore, a quasi-homogeneous t-norm need not be pseudo-homogeneous since $T_D$ is a quasi-homogeneous t-norm but not pseudo-homogeneous (one easily verifies that $T_D$ is a quasi-homogeneous t-norm just by letting the function $\psi :[0, 1] \rightarrow[0, 1]$ be defined by $\psi(x)=\begin{cases}
  0  & x \in [0,1),\\1  & x=1
\end{cases}$
and the function $\varphi:[0, 1] \rightarrow[0, \infty)$ be given by $\varphi(x)=x$, respectively), i.e., the concept of a pseudo-homogeneous t-norm is not a generalization of a quasi-homogeneous t-norm. This motivates us to generalize the concepts of both quasi-homogeneous and pseudo-homogeneous t-norms by introducing the concept of a general pseudo-homogeneous t-norm in which the function $F$ is just a general function without the continuous and boundary conditions and characterize the general pseudo-homogeneous t-norm.

The rest of this article is organized as follows. In Section $2$, we present some elementary concepts and results related to t-norms. In Section $3$, we first introduces the concept of a general pseudo-homogeneous t-norm that generalizes the definitions of both quasi-homogeneous and pseudo-homogeneous t-norms, and then give a representation theorem of a general pseudo-homogeneous t-norm. A conclusion is drawn in Section $4$.

\section{Preliminaries}

In this section, we recall some concepts and results which will be used in the sequel.

\begin{definition}[see\cite{Klement2000}, Definition 1.1]\label{def2.1}
\emph{A function $T:[0,1]^{2}\rightarrow[0,1]$ is a \emph{t-norm} if, for all $x,y\in[0,1]$, it satisfies the following conditions:
\begin{enumerate}
\item [(T1)] $T(x,y)=T(y,x)$.
\item [(T2)] $T(x,T(y,z))=T(T(x,y),z)$.
\item [(T3)] $T(x,y)\leq T(x,z)$ whenever $y\leq z$.
\item [(T4)] $T(x,1)=x$.
\end{enumerate}}
\end{definition}

The following are four basic t-norms $T_{M}$, $T_{P}$, $T_{L}$ and $T_{D}$ given by, respectively:

$T_{M}(x,y)=\min\{x,y\}.$ \ \ \ \ \ \ \ \ \ \ \ \ \ \ \ \ \ $T_{P}(x,y)=xy.$

$T_{L}(x,y)=\max\{x+y-1,0\}.$  \ \ \ \ \ \  $T_{D}(x,y)=\left\{\begin{array}{ll}
0 & {\mbox{ if  } x,y \in[0, 1)},\\
\min\{x,y\} & {\mbox{ otherwise}. }
\end{array}
\right.$

\begin{definition}[see\cite{Xie2016}, Definition 3.1]\label{def2.2}
\emph{A t-norm $T:[0,1]^2\rightarrow[0,1]$ is said to be} pseudo-homogeneous \emph{if it satisfies
\begin{equation*}\label{eq2.0}
T(\lambda x, \lambda y)=F(\lambda, T(x,y))
\end{equation*}
for all $x, y, \lambda \in[0,1]$, where $F:[0,1]^{2}\rightarrow[0,1]$ is a continuous and increasing function with $F(x,1)=0\Leftrightarrow x=0$.}
\end{definition}

\begin{lemma}[see\cite{Klement2000}, Proposition 1.9]\label{lemma2.2} \ \par
\begin{enumerate}
\item [(\romannumeral 1)] The only t-norm $T$ satisfying $T(x,x)=x$ for all $x\in[0,1]$ is the minimum $T_{M}$.
\item [(\romannumeral 2)] The only t-norm $T$ satisfying $T(x,x)=0$ for all $x\in[0,1)$ is the drastic product $T_{D}$.
\end{enumerate}
\end{lemma}

\begin{lemma}[see\cite{Xie2016}, Theorem 3.2]\label{lemma2.1}
$T:[0,1]^2\rightarrow[0,1]$ is a pseudo-homogeneous t-norm satisfying $T(x,x)<x$ for all $x\in(0,1)$ if and only if one of the following items holds:
\begin{enumerate}
\item [(\romannumeral 1)] $$T(x,y)=xy=T_P(x,y)$$ and $$F(x,y)=x^2y.$$
\item [(\romannumeral 2)] For some $\beta <0$, $$T(x,y)=\begin{cases}
  (x^{\beta}+y^{\beta}-1)^{\frac{1}{\beta}} & \mbox{\emph{if} }(x,y)\in(0,1]^2,\\ 0 &\mbox{\emph{otherwise}}\\
  \end{cases}$$
and $$F(x,y)=\begin{cases}
  (x^{\beta}+(xy)^{\beta}-1)^{\frac{1}{\beta}} & \mbox{\emph{if} }(x,y)\in(0,1]^2,\\ 0  &\mbox{\emph{otherwise}}.\\
  \end{cases}$$
\end{enumerate}
\end{lemma}

Let $N$ be a positive integer and $R$ be the set of all real numbers. A function $f:R^N\rightarrow R$ is called additive if and only if it satisfies Cauchy's functional equation
$$f(x+y)=f(x)+f(y)$$
for all $x,y\in R^N$.
\begin{lemma}[see\cite{Kuczma2009}, Theorem 5.5.2]\label{lemma2.3}
If $f:R^N\rightarrow R$ is a continuous additive function, then there exists a $c\in R^N$ such that
$$f(x)=cx$$
where $cx=\sum_{i=1}^{N}c_ix_i$ $(c=(c_1, \cdots, c_N), x=(x_1, \cdots, x_N))$ denotes the scalar product.
\end{lemma}

\begin{lemma}[see\cite{Kuczma2009}, Theorem 13.1.6]\label{lemma2.4}
Let $D$ be one of the sets $(0,1), [0,1), (-1,1)$, $(-1,0)\cup(0,1), (1, \infty), (0, \infty), [0,\infty), (-\infty,0)\cup (0,\infty)$ or $R$. A function $f:D\rightarrow R$ is a continuous solution of the multiplicative Cauchy equation $f(xy)=f(x)f(y)$ if and only if either $f=0$, or $f=1$, or for any $x\in D$, $f$ has one of the following forms:
$f(x)=|x|^c,$
$f(x)=|x|^c sgn(x)$
with a certain $c\in R$. If $0\in D$, then $c>0$.
\end{lemma}

\begin{lemma}[see\cite{Klement2000}, Proposition 2.12]\label{lemma3.1}
For a t-norm $T$ the following are equivalent:
\begin{enumerate}
\item [(\romannumeral 1)] $T$ is Archimedean.
\item [(\romannumeral 2)] $T$ satisfies the limit property (LP: for all $x\in(0,1), \lim\limits_{n \rightarrow\infty}x^{(n)}_{T}=0$).
\item [(\romannumeral 3)] $T$ has only trivial idempotent elements and, whenever
$$\lim\limits_{x \searrow x_0}T(x,x)=x_0$$
for some $x_0\in(0,1)$, there exists a $y_0\in(x_0,1)$ such that $T(y_0,y_0)=x_0$.
\end{enumerate}
\end{lemma}

\begin{lemma}[see\cite{Klement2000}, Proposition 2.15]\label{lemma3.2}
For an arbitrary $t$-norm $T$ we have:
\begin{enumerate}
\item [(\romannumeral 1)] If $T$ is right-continuous and has only trivial idempotent elements then it is Archimedean.
\item [(\romannumeral 2)] If $T$ is right-continuous and satisfies the conditional cancellation law $(CCL)$ then it is Archimedean.
\item [(\romannumeral 3)] If $\lim\limits_{x \searrow x_0}T(x,x)<x_0$ for all $x_0\in (0,1)$ then $T$ is Archimedean.
\item [(\romannumeral 4)] If $T$ is strict then it is Archimedean.
\item [(\romannumeral 5)] If each $x\in (0,1)$ is a nilpotent element of $T$ then $T$ is Archimedean.
\end{enumerate}
\end{lemma}

\section{General pseudo-homogeneous t-norms}

In this section, we first introduce the concept of a general pseudo-homogeneous t-norm, and then devote to characterizing general pseudo-homogeneous t-norms.

\begin{definition}\label{def2.4}
\emph{A t-norm $T:[0,1]^2\rightarrow[0,1]$ is said to be} general pseudo-homogeneous \emph{for a binary function $F:[0,1]^{2}\rightarrow[0,1]$ if for all $x, y, \lambda \in[0, 1]$ it satisfies
\begin{equation}\label{eq2.2}
T(\lambda x, \lambda y)=F(\lambda, T(x,y)).
\end{equation}}
\end{definition}

Clearly, a quasi-homogeneous t-norm is general pseudo-homogeneous, and Definition \ref{def2.2} is a special case of Definition \ref{def2.4}. However, a general pseudo-homogeneous t-norm may be neither quasi-homogeneous nor pseudo-homogeneous (see Figure 1 for their relationship). For example, the t-norm $T_L$ is general pseudo-homogeneous for the function $F:[0,1]^{2}\rightarrow[0,1]$ given by
$$F(x,y)=\max\{x+xy-1,0\}.$$
But it is neither pseudo-homogeneous by Definition \ref{def2.2} nor quasi-homogeneous.
\begin{center}
\includegraphics[height=0.3\textwidth]{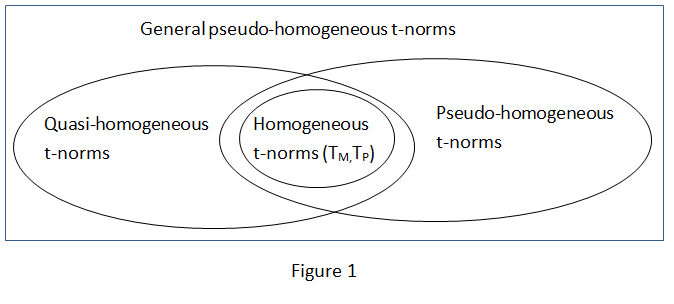}
\end{center}

The following three propositions reveal the relationship between a general pseudo-homogeneous t-norm $T$ and its related function $F$.
\begin{proposition}\label{proposition3.1}
Let $T:[0,1]^2 \rightarrow [0,1]$ be a general pseudo-homogeneous t-norm for a binary function $F:[0,1]^2 \rightarrow[0,1]$. Then the t-norm $T$ and the binary function $F$ are uniquely determined each other.
\end{proposition}
\begin{proof}
Suppose that there exist two functions $F_1:[0,1]^2 \rightarrow [0,1]$ and $F_2:[0,1]^2 \rightarrow [0,1]$ such that
\begin{equation}\label{eq3.1}
T(\lambda x, \lambda y)=F_1(\lambda, T(x, y))
\end{equation}
and
\begin{equation}\label{eq3.2}
T(\lambda x,\lambda y)=F_2(\lambda, T(x,y))
\end{equation}
for all $x, y\in[0,1]$, respectively.
Then for all $x, y\in [0,1]$
\begin{equation*}\label{eq3.2001}
\begin{split}
F_1(x,y)&=F_1(x,T(1,y))\\
&=T(x,xy) \  \ by\  Eq.(\ref{eq3.1})\\
&=F_2(x,T(1,y)) \  \ by\  Eq.(\ref{eq3.2})\\
&=F_2(x,y).
\end{split}
\end{equation*}
Therefore, for the general pseudo-homogeneous t-norm $T$ the binary function $F$ is unique.

On the other hand, suppose that there exist two t-norms $T_1:[0,1]^2 \rightarrow [0,1]$ and $T_2:[0,1]^2 \rightarrow [0,1]$ such that
 \begin{equation}\label{eq3.3}
 T_1(\lambda x, \lambda y)=F(\lambda, T_1(x,y))
 \end{equation}
 and
 \begin{equation}\label{eq3.4}
 T_2(\lambda x, \lambda y)=F(\lambda, T_2(x, y))
 \end{equation}
for all $ x, y, \lambda\in[0,1]$, respectively. Then for all $x,y\in[0,1]$ with $x\geq y$ and $x\neq 0$ we have
 \begin{equation*}\label{equ4001}
  \begin{split}
  T_1(x,y)&=T_1(x\cdot1,x\cdot \frac {y}{x})\\
  &=F(x,T_1(1,\frac {y}{x})) \  \ by\  Eq.(\ref{eq3.3})\\
  &=F(x,\frac {y}{x})\\
  &=F(x,T_2(1,\frac {y}{x})) \\
  &=T_2(x,y)\  \ by \  Eq.(\ref{eq3.4}),
 \end{split}
\end{equation*}
i.e., $T_{1}(x,y)=T_{2}(x,y)$. Obviously, $T_1(0,0)=T_2(0,0)$. From the commutativity of $T$ we can get $T_1(x,y)=T_2(x,y)$ for all $x, y\in[0,1]$ with $x\leq y$ and $y\neq 0$. Consequently, for the binary function $F$ the general pseudo-homogeneous t-norm $T$ is unique.
\end{proof}

\begin{proposition}\label{proposition3.2}
Let $T:[0,1]^2 \rightarrow [0,1]$ be a general pseudo-homogeneous t-norm for a binary function $F:[0,1]^2 \rightarrow[0,1]$. Then $T$ is continuous if and only if $F$ is continuous.
\end{proposition}
\begin{proof}
Suppose that the t-norm $T$ is continuous. From Definition \ref{def2.4}, we have
\begin{equation*}\label{equ4002}
F(x,y)=F(x,T(1,y))=T(x,xy)
\end{equation*}
for all $x, y \in[0,1]$. Thus, $F$ is continuous.

On the other hand, if the binary function $F$ is continuous, then we have
\begin{equation*}\label{equ4003}
\begin{split}
T(x,y)&=T(x\cdot1, x\cdot \frac{y}{x})\\
&=F(x, T(1,\frac{y}{x}))\\
&=F(x,\frac{y}{x})
\end{split}
\end{equation*}
for all $x, y\in[0,1]$ with $x\geq y$ and $x\neq 0$.
Similarly, we can get $T(x,y)=F(y,\frac{x}{y})$ for all $x, y\in[0,1]$ with $x\leq y$ and $y\neq0$.
Thus, for all $x, y\in[0,1]$ with $(x,y)\neq (0,0)$, we have
\begin{equation*}\label{equ4004}
T(x, y)=F(\max\{x,y\},\frac{\min\{x,y\}}{\max\{x,y\}}).
\end{equation*}
So, to show the continuity of $T$, it is enough to prove that $T$ is continuous at the point $(0,0)$. In fact,
\begin{equation*}\label{equ4005}
0\leq\displaystyle\lim_{(x,y)\searrow(0,0)}T(x, y)\leq\displaystyle\lim_{(x,y)\searrow(0,0)}\min\{x, y\}=0.
\end{equation*}
This shows that
\begin{equation*}\label{equ4006}
\displaystyle\lim_{(x,y)\searrow(0,0)}T(x,y)=0=T(0,0).
\end{equation*}
Therefore, $T$ is continuous at the point $(0,0)$.
\end{proof}

\begin{proposition}\label{prop4.1}
Let $T:[0,1]^2 \rightarrow [0,1]$ be a general pseudo-homogeneous t-norm for a binary function $F:[0,1]^2 \rightarrow[0,1]$. Then the following statements are equivalent:
\begin{enumerate}
\item [(\romannumeral 1)] $T=T_M$.
\item [(\romannumeral 2)] $F(x,y)=xy$ for all $x,y\in[0,1]$.
\item [(\romannumeral 3)] $F$ is commutative.
\item [(\romannumeral 4)] $F(x,1)=x$ for all $x\in[0,1]$.
\end{enumerate}
\end{proposition}
\begin{proof} $(\romannumeral 1) \Rightarrow (\romannumeral 2)$ Let $T=T_M$. From Eq.(\ref{eq2.2}) we have
$F(x,y)=F(x,T(1,y))=T(x,xy)=\min\{x,xy\}=xy$ for all $x,y\in[0,1]$.

$(\romannumeral 2) \Rightarrow (\romannumeral 3)$ Trivially, $(\romannumeral 2)$ implies $(\romannumeral 3)$.

$(\romannumeral 3) \Rightarrow (\romannumeral 4)$ $F(x,1)=F(1,x)=F(1,T(1,x))=T(1,x)=x$ for all $x\in[0,1]$ since $F$ is commutative.

$(\romannumeral 4) \Rightarrow (\romannumeral 1)$ From Eq.(\ref{eq2.2}) we have $F(x,y)=F(x,T(1,y))=T(x,xy)$ for all $x,y\in[0,1]$. Then
$x=F(x,1)=T(x,x), x\in[0,1]$. It follows immediately from Lemma \ref{lemma2.2} (i) that $T=T_M$.
\end{proof}

A function $\delta_{T}: [0, 1] \rightarrow[0, 1]$ with
$\delta_{T}(x) = T(x, x)$ is called a diagonal function of a t-norm $T$, in symbols $\delta_T$. Then we have the following proposition.

\begin{proposition}\label{lema3.1}
Let $T:[0,1]^2\rightarrow[0,1]$ be a general pseudo-homogeneous t-norm. Then the following items are equivalent:
\begin{enumerate}
\item [(\romannumeral 1)] The diagonal $\delta_T$ is left-continuous on $(0, 1)$.
\item [(\romannumeral 2)] The diagonal $\delta_T$ is continuous on $(0, 1)$.
\end{enumerate}
\end{proposition}
\begin{proof} Trivially, (\romannumeral 2) implies $(\romannumeral 1)$. In order to show the converse, we only need to show that $\delta_T$ must be right-continuous when $T$ is general pseudo-homogeneous. In fact, for any $x_0\in (0,1)$ and $x_0\leq y,z<1$ with $x_0=yz$, we have
\begin{equation*}\label{eq0002}
T(x_0,x_0)=T(yz,yz)=T(y,yT(z,z)).
\end{equation*}
Then
$$T(x_0,x_0)=\lim\limits_{z\searrow x_0(x_0=yz,y \nearrow 1)}T(yz,yz)=\lim\limits_{z\searrow x_0(x_0=yz,y \nearrow 1)}T(y,yT(z,z))=T(x^{+}_0,x^{+}_0),$$ i.e., $\delta_T$ is right-continuous.
\end{proof}

Denoted by $A\setminus B=\{x\in A\mid x\notin B\}$ for two sets $A$ and $B$. Then we have the following proposition.
\begin{proposition}\label{Proposition4.4}
Let $T:[0,1]^2\rightarrow[0,1]$ be a general pseudo-homogeneous t-norm for a binary function $F:[0,1]^2 \rightarrow[0,1]$. If $\lim\limits_{x \searrow x_0}T(x,x)=x_0$ for an $x_0\in (0,1)$ then the following items hold.
\begin{enumerate}
\item [(\romannumeral 1)] $\delta_T(x)=x$ for all $x\in (x_0,1)$.
\item [(\romannumeral 2)] $T(x,y)=\min\{x,y\}$ for all $x, y\in (x_0,1]$.
\item [(\romannumeral 3)] One of the following two cases holds:
\begin{enumerate}
\item [$(1^\star)$] $T=T_M$.
\item [$(2^\star)$] $T(x_0,x_0)=x_0$ and there exists $0< c <1$ such that
$$T(x,y)=\begin{cases}
0  & \emph{if }(x,y)\in(0,1)^2\setminus[c,1)^2,\\  \min\{x,y\}  & \emph{otherwise}.
\end{cases}$$
  \end{enumerate}
\end{enumerate}
\end{proposition}
\begin{proof}
$(\romannumeral 1)$ Suppose that there exists a $y\in (x_0,1)$ such that $T(y,y)=y_0<y$. Let $\lambda =\frac{x}{y}$. Then
\begin{equation*}\label{equ4011}
  \begin{split}
  x_0&=\lim\limits_{x \searrow x_0}T(x,x)\\
  &=\lim\limits_{\lambda\searrow \frac{x_0}{y}}T(\lambda y,\lambda y) \\
  &=\lim\limits_{\lambda\searrow \frac{x_0}{y}}T(\lambda ,\lambda T(y,y)) \\
  &=\lim\limits_{\lambda\searrow \frac{x_0}{y}}T(\lambda ,\lambda y_0) \\
  &\leq\lim\limits_{\lambda\searrow \frac{x_0}{y}}\lambda y_0 \\
  &=\lim\limits_{x \searrow x_0}\frac{x}{y}\cdot y_0 \\
  &=\frac{y_0}{y}\cdot x_0  \\
  &<x_0,
 \end{split}
\end{equation*}
a contradiction. So $T(x,x)=x$ for all $x\in (x_0,1)$, i.e., $\delta_T(x)=x$ for all $x\in (x_0,1)$.

$(\romannumeral 2)$ Obviously, from $(\romannumeral 1)$ and Definition \ref{def2.1}, we get $T(x,y)=y$ for all $x_0<y\leq x\leq 1$ and
$T(x,y)=x$ for all $x_0<x\leq y \leq 1$. Then for all $(x,y)\in (x_0,1]^2$ we have
\begin{equation*}\label{eq1}
T(x,y)=\min\{x,y\}.
\end{equation*}

$(\romannumeral 3)$ From Definition \ref{def2.4} and $(\romannumeral 2)$, for all $x_0\leq y,z<1$ with $x_0=yz$, we have
\begin{equation*}\label{eq2}
T(x_0,x_0)=T(yz,yz)=T(y,yT(z,z)).
\end{equation*}
Then
\begin{equation*}\label{eq3}
 \begin{split}
T(x_0,x_0)&=\lim\limits_{z\searrow x_0(x_0=yz,y \nearrow 1)}T(yz,yz)\\
&=\lim\limits_{z\searrow x_0(x_0=yz,y \nearrow 1)}T(y,yT(z,z))\\
&=T(x^{+}_0,x^{+}_0)\\
&=\lim\limits_{x\searrow x_0}T(x,x)\\
&=x_0.
 \end{split}
\end{equation*}
Thus either $\lim\limits_{y \searrow x}T(y,y)=x $ for all $x\in (0,1)$
or
there exists a $c$ with $0<c\leq x_0$ such that
$T(c,c)=c, \lim\limits_{y \searrow x}T(y,y)=x \mbox{ for all } x\in (c,1)\mbox{ and }\lim\limits_{z \searrow x}T(z,z)< x \mbox{ for all } x\in (0,c)$.

Case $(1^\star)$. If $\lim\limits_{y \searrow x}T(y,y)=x $ for all $x\in (0,1)$, then from (i) we have $T(x,x)=x$ for all $x\in (0,1)$. It follows immediately from Lemma \ref{lemma2.2} (i) that $T=T_M$.

Case $(2^\star)$. If there exists a $c$ with $0<c\leq x_0$ such that
$T(c,c)=c, \lim\limits_{y \searrow x}T(y,y)=x \mbox{ for all } x\in (c,1)\mbox{ and }\lim\limits_{z \searrow x}T(z,z)<x \mbox{ for all } x\in (0,c)$, then we can take $x_c$ arbitrarily such that $x_c$ is sufficiently close to $c$ with $x_c<c$ and $T(x_c,x_c)<x_c$.
Let $T(x_c,x_c)=a$. Then from Definition \ref{def2.4} and $(\romannumeral 2)$, for all $c\leq y,z<1$ with $x_c=yz$, we have
\begin{equation*}\label{eq00200}
a=T(x_c,x_c)=T(yz,yz)=T(y,yT(z,z))=T(y,yz)=T(y,x_c)=T(x_c,y).
\end{equation*}
Hence $T(x_c,T(y,z))=T(x_c, \min\{y,z\})$ and $T(x_c,T(y,z))=T(T(x_c,y),z)=T(a,z)$ for all $c\leq y, z<1$ and $x_c=yz$.
Thus $T(z,a)=T(a,z)=a, c\leq z<1.$ Then, for all $c\leq z<1$ and $a\leq x\leq x_c$, we have
\begin{equation}\label{eq7}
T(z,x)=T(x,z)=a.
\end{equation}
Take a number $b$ such that $a<b<x_c$ with $b=y_1 z_1$ and $c\leq y_1, z_1<1$. Thus
\begin{equation*}\label{eq8}
T(b,b)=T(z_1 y_1,z_1 y_1)=T(z_1,z_1 T(y_1,y_1))=T(z_1,z_1 y_1)=T(z_1,b).
\end{equation*}
Then $\delta_T(b)=T(b,b)=a$. Hence for all $x, y\in [b,x_c]$,
\begin{equation*}\label{equ4014}
T(x,y)=T(x_c,x_c)=T(b,b)=\delta_T(b)=a.
\end{equation*}
 Let $\lambda_1=\frac{b}{x_c}$. Then for all $x, y\in [b,x_c]$,
\begin{equation*}\label{equ4015}
T(\lambda_1 x,\lambda_1 y)=F(\lambda_1, T(x, y)) =F(\lambda_1, \delta_T(b))
\end{equation*}
and $\lambda_1 x, \lambda_1 y\in [\frac {b^2}{x_c}, b]$.
Thus $T(x,y)=F(\lambda_1, \delta_T(b))$ for all $x, y\in [\frac {b^2}{x_c}, b]$.
Then $$\delta_T(\frac {b^2}{x_c})= T(\frac {b^2}{x_c},\frac {b^2}{x_c})=T(x_c,x_c)=\delta_T(x_c)=a$$
and
\begin{equation*}\label{equ4016}
T(x,y)=\delta_T(\frac {b^2}{x_c})=\delta_T(b)=a, x,y\in [\frac {b^2}{x_c},b].
\end{equation*}
Similarly, let $\lambda_2=\frac {b^2}{x_c^2}$. Then for all $x, y\in [\frac {b^2}{x_c}, x_c]$,
\begin{equation*}\label{equ4017}
T(\lambda_2 x,\lambda_2 y)=F(\lambda_2, T(x,y))=F(\lambda_2,\delta_T(b))
\end{equation*}
and $\lambda_2 x, \lambda_2 y\in [\frac {b^4}{x_c^3},  \frac {b^2}{x_c}]$.
Thus $T(x,y)=F(\lambda_2, \delta_T(b))$ for all $x, y\in [\frac {b^4}{x_c^3}, \frac {b^2}{x_c}]$.
Then
$$\delta_T(\frac{b^4}{x_c^3})=T(\frac{b^4}{x_c^3},\frac {b^4}{x_c^3})=T(\frac{b^2}{x_c},\frac {b^2}{x_c})=\delta_T(\frac{b^2}{x_c})=\delta_T(b)$$
and
$$T(x,y)=\delta_T(\frac {b^4}{x_c^3})=\delta_T(b), x, y\in [\frac {b^4}{x_c^3}, x_c].$$
Continuing the above discussion, for all $\lambda_n=\frac {b^{2^{(n-1)}}}{x_c^{2^{(n-1)}}}$ ($n=1, 2, 3, \cdots,$)
we obtain a closed interval sequence $([\frac {b^{2^{n}}}{x_c^{2^{n}-1}}, b])_{n=1, 2, 3, \cdots,}$ which satisfies
$$\delta_T(\frac{b^{2^{n}}}{x_c^{2^{n}-1}})=\delta_T(b) \mbox{ and }\lim\limits_{n\rightarrow\infty} \frac{b^{2^{n}}}{x_c^{2^{n}-1}}=0.$$
Then
\begin{equation*}\label{equ4018}
\begin{split}
0&=T(0,0)\\
&=T(\lim \limits_{n\rightarrow\infty} \frac{b^{2^{n}}}{x_c^{2^{n}-1}}, \lim \limits_{n\rightarrow\infty} \frac{b^{2^{n}}}{x_c^{2^{n}-1}}) \\
&=\lim \limits_{n\rightarrow\infty}T(\frac{b^{2^{n}}}{x_c^{2^{n}-1}},\frac{b^{2^{n}}}{x_c^{2^{n}-1}})\\
&=\lim \limits_{n\rightarrow\infty}\delta_T(\frac{b^{2^{n}}}{x_c^{2^{n}-1}})\\
&=\lim \limits_{n\rightarrow\infty}\delta_T(b)\\
&=\delta_T(b).
\end{split}
\end{equation*}
Therefore, $a=\delta_T(b)=0$, i.e., $T(x_c,x_c)=0$. So that
\begin{equation}\label{eq10}
T(x,y)=0
\end{equation}
 for all $0\leq x \leq x_c, 0\leq y \leq x_c$. Then from $(\romannumeral 2)$ and Eqs. (\ref {eq7}) and (\ref {eq10}), we have
$$T(x,y)=\begin{cases}
0 & \mbox{if }(x,y)\in(0,1)^2\setminus[c,1)^2,\\  \min\{x,y\}  & \mbox{otherwise}.
  \end{cases}$$
\end{proof}

In order to present a representation theorem of a general pseudo-homogeneous t-norm, the following four propositions are needed also.
\begin{proposition}\label{proposition1}
Let $T:[0,1]^2\rightarrow[0,1]$ be a general pseudo-homogeneous t-norm for a binary function $F:[0,1]^2 \rightarrow[0,1]$. Then $\lim\limits_{x \nearrow 1}T(x,x)=1 $ or $\lim\limits_{x \nearrow 1}T(x,x)=0$.
\end{proposition}
\begin{proof}
Obviously, for a general pseudo-homogeneous continuous t-norm $T$, $\lim\limits_{x \nearrow 1}T(x,x)=1$ holds. Next, we show that $\lim\limits_{x \nearrow 1}T(x,x)=0$ if $\lim\limits_{x \nearrow 1}T(x,x)<1$. Suppose that $\lim\limits_{x \nearrow 1}T(x,x)=a$ with $a<1$. Then from the proof of Proposition \ref{Proposition4.4} we have $\lim\limits_{x \searrow x_0}T(x,x)<x_0$ for all $x_0\in(0,1)$. Therefore, by Lemma \ref{lemma3.2} (iii), $T$ is Archimedean. Furthermore, from Lemma \ref{lemma3.1}
\begin{equation}\label{eq100}
\lim\limits_{n \rightarrow\infty}x^{(n)}_{T}=0
\end{equation}
for all $x\in(0,1)$. On the other hand, for any $x\in (0,1)$, if $x=yz$ with $y,z\in [x,1]$ then
\begin{equation*}\label{eq101}
T(x,x)=T(yz,yz)=T(z,zT(y,y)).
\end{equation*}
Thus
\begin{equation*}\label{eq102}
T(x,x)=\lim\limits_{y \nearrow 1(x=yz, z\searrow x)}T(yz,yz)=\lim\limits_{y \nearrow 1(x=yz,z\searrow x)}T(z,zT(y,y))=T(x,xa),
\end{equation*}
i.e., $T(x,x)=T(x,xa)$ for all $x\in (0,1)$, which results in
\begin{equation*}\label{eq104}
\begin{split}
T(T(x,x),x)&=T(x,T(x,x))\\
&=T(x,T(x,xa))\\
&=T(T(x,x),xa)\\
&=T(T(x,xa),xa)\\
&=T(x,T(xa,xa))\\
&=T(x,T(x,xT(a,a)))\\
&=T(T(x,x),xT(a,a)).
\end{split}
\end{equation*}
Consequently, $T(T(x,x),x)=T(T(x,x),T(x,x))$, i.e., $x^{(3)}_{T}=x^{(4)}_{T}$ because $xT(a,a)<T(x,x)< x$ when $a\leq x$.
Thus $x^{(n)}_{T}=x^{(3)}_{T}$ for all $n\geq 4$ and $x\in [a,1]$, which together with Eqs. (\ref{eq100}) implies that $x^{(3)}_{T}=x^{(n)}_{T}=0$ for all $x\in [a,1]$. So that for any $y\in (0,a)$ there exits an $x\in [a,1]$ such that $y\leq T(x,x)$ since $\lim\limits_{x \nearrow 1}T(x,x)=a$. Thus $T(x,y)\leq T(x,T(x,x))=0$, implying $T(x,y)=0$. Then $T(z,y)\leq T(z,T(z,z))=0$, i.e., $T(z,y)=0$ for all $x\leq z<1$ and $T(t,y)\leq T(x,y)=0$, i.e., $T(t,y)=0$ for all $0<t<x$. Applying Definition \ref{def2.1} (T1), we obtain $T(x,y)=0$ for all $(x,y)\in(0,1)^2\setminus[a,1)^2$.
Therefore, $T(x,x)=T(x,ax)=0$ for all $x\in[a,1)$ since $ax<a$, which means $\lim\limits_{x \nearrow 1}T(x,x)=0$.
\end{proof}

\begin{proposition}\label{Proposition4.6}
Let $T:[0,1]^2\rightarrow[0,1]$ be a general pseudo-homogeneous t-norm for a binary function $F:[0,1]^2 \rightarrow[0,1]$. If $\lim\limits_{x \searrow x_0}T(x,x)< x_0$ for all $x_0\in (0,1)$, then the following items holds:
\begin{enumerate}
\item [(\romannumeral 1)] The diagonal $\delta_T$ is a continuous increasing function on $(0,1)$.
\item [(\romannumeral 2)] $T=T_D$ or $T$ is a continuous Archimedean t-norm.
\end{enumerate}
\end{proposition}
\begin{proof}
From Lemma \ref{lemma3.2} we know that the t-norm $T$ is Archimedean.

Concerning $(\romannumeral 1)$, trivially, the diagonal $\delta_T$ is increasing on $(0,1)$. Concerning the continuity of $\delta_T$, it suffices to show that $\delta_T$ is left-continuous because of Lemma \ref{lema3.1}. From Proposition \ref{proposition1}, we distinguish by two cases as follows.

Case 1. If $\lim\limits_{x \nearrow 1}T(x,x)=0$, then $\delta_T(x)=T(x,x)=0$ for all $x\in (0,1)$, i.e., $\delta_T$ is a continuous increasing function on $(0,1)$.\\

Case 2. Fix an arbitrary strictly increasing sequence $(z_n)_{n\in N} \in [0,1]$ with $\lim\limits_{n \rightarrow \infty}z_n=1.$ If $\lim\limits_{x \nearrow 1}T(x,x)=1$, then $\lim\limits_{n \rightarrow \infty}(z_n)_T^{(2)}=1$. Assume that $T$ is general pseudo-homogeneous and Archimedean but its diagonal $\delta_T(x)$ is not left-continuous in some point $x_0\in (0,1)$. Then for each $n\in N$ there exist numbers $k_n\in N$ such that
$$(z_n)_T^{(k_n)}\leq x_0<(z_n)_T^{(k_n-1)}$$
implying that for all $n\in N$
$$(z_n)_T^{(2k_n)}\leq T(x_0^-, x_0^-)<T(x_0, x_0)\leq (z_n)_T^{(2k_n-2)}$$
since $T$ is Archimedean. $\lim\limits_{n \rightarrow \infty}(z_n)_T^{(2)}=1$ and $x_0<(z_n)_T^{(k_n-1)}$ yield that
there exists a number $n\in N$ such that $(z_n)_T^{(2)}\cdot(z_n)_T^{(k_n-1)}\geq x_0$
\begin{equation*}\label{eq110}
\begin{split}
(z_n)_T^{(2k_n)}&=T((z_n)_T^{(2)},(z_n)_T^{(2k_n-2)})\\
&\geq T((z_n)_T^{(2)},(z_n)_T^{(2)}\cdot(z_n)_T^{(2k_n-2)})\\
&=T((z_n)_T^{(2)}\cdot(z_n)_T^{(k_n-1)},(z_n)_T^{(2)}\cdot(z_n)_T^{(k_n-1)})\\
&\geq T(x_0, x_0)\\
&>T(x_0^-, x_0^-)\\
&\geq(z_n)_T^{(2k_n)}
\end{split}
\end{equation*}
which is a contradiction.

Concerning $(\romannumeral 2)$, trivially, if $\lim\limits_{x \nearrow 1}T(x,x)=0$ then $T=T_D$. If $\lim\limits_{x \nearrow 1}T(x,x)=1$, then $\delta_T$ is continuous. Next, we show that the pseudo-inverse $\delta_T^{(-1)}$ of $\delta_T$ is continuous. Due to Remark 3.4 (ii) in \cite{Klement2000}, we only need to show that $\delta_T$ is strictly monotone on the set $\delta_T^{(-1)}([0,1))$. In fact, if there exists a subinterval $[a,b]$ of $[0,1]$ such that $\delta_T(a)=\delta_T(b)$, i.e., $\delta_T$ is a constant function on $[a,b]$, then we have $\delta_T(b)=\delta_T(a)=\delta_T(0)=0$ by simply imitating the proof of Case $(2^\star)$ in Proposition \ref{Proposition4.4} (iii). Thus $\delta_T$ can be given as
$$\delta_T(x)= \left\{\begin{array}{ll}
0 & {\mbox{\scriptsize\normalsize  if } \   0 \leq x \leq c,}\\
f_{1}(x) & {\mbox{\scriptsize\normalsize if }\   c \leq x \leq 1}
\end{array}
\right.$$
where $c$ is a constant with $0\leq c<1$ and $f_{1}$ is strict. This follows that $\delta_T^{(-1)}([0,1))=[c,1)$ and $\delta_T$ is strictly monotone on the set $[c,1)$. So that the pseudo-inverse $\delta_T^{(-1)}$ of $\delta_T$ is continuous.
Because
\begin{equation*}\label{eq117}
\begin{split}
F(x,y)&=F(x,\delta_T(\delta_T^{(-1)}(y)))\\
&=F(x,T(\delta_T^{(-1)}(y),\delta_T^{(-1)}(y)))\\
&=T(x\delta_T^{(-1)}(y),x\delta_T^{(-1)}(y))\\
&=\delta_T(x\delta_T^{(-1)}(y))
\end{split}
\end{equation*}
 for all $x,y\in [0,1]$, implying $F$ is continuous. A direct consequence of Proposition \ref{proposition3.2} is that $T$ is a continuous Archimedean t-norm.
\end{proof}

\begin{proposition}\label{prop4.2}
If a continuous t-norm $T$ is general pseudo-homogeneous for a binary function $F:[0,1]^2 \rightarrow[0,1]$ then $T$ cannot be a non-trivial ordinal sum of continuous Archimedean t-norms.
\end{proposition}
\begin{proof}
Suppose that $T=(\langle a_\alpha,e_\alpha,T_\alpha\rangle)_{\alpha\in A}$ and there exists a summand $\langle a_{\alpha_0},e_{\alpha_0},T_{\alpha_0}\rangle$, where $a_{\alpha_0}>0$ and $e_{\alpha_0}\leq1$, or $a_{\alpha_0}\geq0$ and $e_{\alpha_0}<1$, $\alpha_0\in A$ and $T_{\alpha_0}$ is a continuous Archimedean t-norm. Moreover, from Definition \ref{def2.4} we get
\begin{equation*}\label{equ4023}
T(\lambda x,\lambda y)=F(\lambda,T(x,y))
\end{equation*}
for all $x, y, \lambda \in[0,1]$ and $F(x,y)=T(x,xy)$. Next, we investigate the t-norm $T$ by distinguishing two cases.\\
Case 1. $a_{\alpha_0}>0$ and $e_{\alpha_0}\leq1$. Then
$$T(x,a_{\alpha_0})=\min\{x,a_{\alpha_0}\}$$
for all $x\in[0,1]$ and there exists $x, y\in(a_{\alpha_0}, e_{\alpha_0})$ such that $a_{\alpha_0}<y<x<e_{\alpha_0}$ and $T(x,y)<y$. Let $\lambda=\frac{a_{\alpha_0}}{y}$.
Then
\begin{equation*}\label{equ4007}
\begin{split}
T(\lambda x, \lambda y)&=T(\frac{a_{\alpha_0}}{y} x, \frac{a_{\alpha_0}}{y} y)\\
&=T(\frac{a_{\alpha_0}}{y} x, a_{\alpha_0})\\
&=\min\{\frac{x}{y}a_{\alpha_0}, a_{\alpha_0}\}\\
&=a_{\alpha_0}
\end{split}
\end{equation*}
and
\begin{equation*}\label{equ4008}
  \begin{split}
  F(\lambda,T(x, y))&=T(\lambda, \lambda T(x, y))\\
  &\leq \lambda T(x, y) \\
  &=\frac{a_{\alpha_0}}{y}T(x,y),
 \end{split}
\end{equation*}
which imply that $a_{\alpha_0}\leq\frac{a_{\alpha_0}}{y}T(x, y)$, i.e., $T(x, y)\geq y$, contrary to $T(x,y)<y$.

Case 2. $a_{\alpha_0}\geq0$ and $e_{\alpha_0}<1$. Then there exists an $x\in(a_{\alpha_0}, e_{\alpha_0})$ such that $e^2_{\alpha_0}<x<e_{\alpha_0}$. Let $\lambda=\frac{x}{e_{\alpha_0}}$.
Then
\begin{equation*}\label{equ4024}
T(\lambda e_{\alpha_0}, \lambda e_{\alpha_0})=T(x, x)<x
\end{equation*}
and
\begin{equation*}\label{equ4009}
  \begin{split}
  F(\lambda,T(e_{\alpha_0}, e_{\alpha_0}))&=T(\lambda, \lambda T(e_{\alpha_0}, e_{\alpha_0}))\\
  &=T(\frac{x}{e_{\alpha_0}}, \frac{x}{e_{\alpha_0}} e_{\alpha_0}) \\
  &=T(\frac{x}{e_{\alpha_0}}, x)\\
  &=x,
 \end{split}
\end{equation*}
which imply that $x<x$, a contradiction.
\end{proof}

The following proposition can bee seen as a representation theorem of a general pseudo-homogeneous t-norm when it is continuous under some conditions.
\begin{proposition}\label{prop0003.5}
A continuous t-norm $T:[0,1]^2\rightarrow[0,1]$, satisfying $T(x,x)<x$ for all $x\in(0,1)$, is general pseudo-homogeneous for a binary function $F:[0,1]^2 \rightarrow[0,1]$ if and only if one of the following items holds.
\begin{enumerate}
\item [(\romannumeral 1)] For some $\beta >0$,
\begin{equation*}\label{equ4025}
T(x,y)=\begin{cases}
  (\max\{x^{\beta}+y^{\beta}-1,0\})^{\frac{1}{\beta}} & \mbox{\emph{if }}(x,y)\in(0,1]^2,\\0  & \emph{otherwise }
  \end{cases}
\end{equation*}
 and
$$F(x,y)=\begin{cases}
  (\max\{x^{\beta}+(xy)^{\beta}-1,0\})^{\frac{1}{\beta}} & \emph{if }(x,y)\in(0,1]^2,\\0  & \emph{otherwise }.
  \end{cases}$$
\item [(\romannumeral 2)] $$T(x,y)=xy=T_P(x,y)$$ and $$F(x,y)=x^2y.$$
\item [(\romannumeral 3)] For some $\beta<0$, $$T(x,y)=\begin{cases}
  (x^{\beta}+y^{\beta}-1)^{\frac{1}{\beta}} & \emph{if }(x,y)\in(0,1]^2,\\0  & \emph{otherwise }
  \end{cases}$$
and $$F(x,y)=\begin{cases}
  (x^{\beta}+(xy)^{\beta}-1)^{\frac{1}{\beta}} & \emph{if }(x,y)\in(0,1]^2,\\0 & \emph{otherwise }.
  \end{cases}$$
\end{enumerate}
\end{proposition}
\begin{proof}
$(\Rightarrow)$ If $T:[0,1]^2\rightarrow[0,1]$ is a general pseudo-homogeneous continuous t-norm satisfying $T(x,x)<x$ for all $x\in(0,1)$ then, obviously, $T$ is continuous Archimedean. Thus, from Theorem 2.8 of \cite{Klement2000}, $T$ is nilpotent or strict. If $T$ is strict then $F(x,1)=T(x,x)=0\Leftrightarrow x=0$. Thus from Definition \ref{def2.2} and Lemma \ref{lemma2.1} we have $(\romannumeral 2)$ and $(\romannumeral 3)$. If $T$ is nilpotent, then there is a strictly increasing bijection $\phi:[0,1]\rightarrow[0,1]$ such that
\begin{equation*}\label{eq4.6}
T(x,y)=\phi^{-1}(T_{L}(\phi(x),\phi(y)))
\end{equation*}
 for all $(x,y)\in[0,1]^2$, which together with Eq.(\ref{eq2.2}) deduces
$$T_L(\phi(\lambda x), \phi(\lambda y))=T_L(\phi(\lambda),\phi(\lambda \phi^{-1}(T_L(\phi(x),\phi(y))))),$$
i.e.,
$$\max\{\phi(\lambda x)+\phi(\lambda y)-1,0\}=\max\{\phi(\lambda)+\phi(\lambda \phi^{-1}(\max\{\phi(x)+\phi(y)-1,0\}))-1,0\}$$
for all $ x,y,\lambda \in[0,1]$. This follows that $\phi(\lambda x)+\phi(\lambda y)-1>0$ if and only if $\phi(\lambda)+\phi(\lambda \phi^{-1}(\max\{\phi(x)+\phi(y)-1,0\}))-1>0$, and
$$\phi(\lambda x)+\phi(\lambda y)-1=\phi(\lambda)+\phi(\lambda \phi^{-1}(\phi(x)+\phi(y)-1))-1.$$
Thus
\begin{equation}\label{eq4.7}
\phi(\lambda x)+\phi(\lambda y)=\phi(\lambda)+\phi(\lambda \phi^{-1}(\phi(x)+\phi(y)-1)).
\end{equation}
For a given $\lambda$, define a function $g_{\lambda}:[0,1]\rightarrow[0,1]$ by $g_{\lambda}(x)=\phi(\lambda \phi^{-1}(x))$. Let $s=\phi(x)$ and $t=\phi(y)$.
Clearly, $g_{\lambda}$ is a strictly increasing bijection, $g_{\lambda}(s)=\phi(\lambda x)$, $g_{\lambda}(t)=\phi(\lambda y)$ and $g_{\lambda}(1)=\phi(\lambda)$.
Therefore, from Eq.(\ref {eq4.7}), we have
$$g_{\lambda}(s)+g_{\lambda}(t)=g_{\lambda}(1)+g_{\lambda}(s+t-1),$$
i.e.,
$$g_{\lambda}(s+t-1)=g_{\lambda}(s)+g_{\lambda}(t)-g_{\lambda}(1).$$
Define a function $G_{\lambda}:R\rightarrow R$ by $G_{\lambda}(x)=g_{\lambda}(x)$ for all $x\in[0,1]$ and $G_{\lambda}(x+y-1)=G_{\lambda}(x)+G_{\lambda}(y)-G_{\lambda}(1)$ for all $x,y\in R$.
Thus
$$G_{\lambda}(x+y)=G_{\lambda}(x)+G_{\lambda}(y)$$
 for $x,y\in R$, i.e., $G_{\lambda}$ satisfies the additive Cauchy equation. Because of Lemma \ref{lemma2.3}, there exists some constant $c_{\lambda}>0$ such that $G_{\lambda}(x)=c_{\lambda}x$ for all $x\in R$.
Consequently, $g_{\lambda}(x)=G_{\lambda}(x)=c_{\lambda}x$ for all $x\in[0, 1]$.
Then
\begin{equation*}
\phi(\lambda \phi^{-1}(x))=c_{\lambda}x.
\end{equation*}
Note that $c_{\lambda}$ can be rewritten as $c(\lambda)$ since it is related with $\lambda$. So we have
\begin{equation}\label{eq4.8}
\phi(\lambda \phi^{-1}(x))=c(\lambda)x.
\end{equation}
Putting $x=1$ in Eq.(\ref {eq4.8}), we have $\phi(\lambda)=c(\lambda)$, which together with Eqs.(\ref {eq4.8}) implies $\phi(\lambda x)=\phi(\lambda)\phi(x)$ for all $x,\lambda\in[0,1]$. This means that $\phi$ satisfies the multiplicative Cauchy equation in [0,1). According to Lemma \ref{lemma2.4}, there exists some constant $\beta>0$ such that $\phi(x)=x^\beta$ for all $x\in[0, 1)$. Obviously, $\phi(x)=x^\beta$ holds if $x=1$. Consequently, $\phi(x)=x^\beta$ for all $x\in[0, 1]$ with the constant $\beta>0$. Therefore,
for some $\beta >0$, $$T(x,y)=\begin{cases}
  (\max\{x^{\beta}+y^{\beta}-1,0\})^{\frac{1}{\beta}} & \mbox{if }(x,y)\in(0,1]^2,\\0  &\mbox{otherwise }
  \end{cases}$$
 and $$F(x,y)=T(x,xy)=\begin{cases}
  (\max\{x^{\beta}+(xy)^{\beta}-1,0\})^{\frac{1}{\beta}} & \mbox{if }(x,y)\in(0,1]^2,\\0  &\mbox{otherwise }.
  \end{cases}$$
$(\Leftarrow)$ It is only a matter of a simple computation.
\end{proof}

As an immediate consequence of Propositions \ref{prop4.1}, \ref{Proposition4.4}, \ref{Proposition4.6}, \ref{prop4.2} and \ref{prop0003.5}, we finally deduce a representation theorem of a general pseudo-homogeneous t-norm as follows.
\begin{theorem}\label{th0003.1}
A t-norm $T:[0,1]^2\rightarrow[0,1]$ is general pseudo-homogeneous for a binary function $F:[0,1]^2 \rightarrow[0,1]$ if and only if one of the following items holds.
\begin{enumerate}
\item [(\romannumeral 1)]
\begin{equation*}\label{eq122}
T(x,y)=\min\{x,y\}=T_M(x,y)
\end{equation*}
and
\begin{equation*}\label{eq123}
F(x,y)=xy.
\end{equation*}
\item [(\romannumeral 2)] For some $\beta >0$,
\begin{equation*}\label{eq121}
T(x,y)=\begin{cases}
  (\max\{x^{\beta}+y^{\beta}-1,0\})^{\frac{1}{\beta}} & \emph{if }(x,y)\in(0,1]^2,\\ 0 & \emph{otherwise}
  \end{cases}
\end{equation*}
and
\begin{equation*}\label{eq119}
F(x,y)=\begin{cases}
  (\max\{x^{\beta}+(xy)^{\beta}-1,0\})^{\frac{1}{\beta}} & \emph{if }(x,y)\in(0,1]^2,\\ 0 & \emph{otherwise}.
\end{cases}
\end{equation*}
\item [(\romannumeral 3)] $$T(x,y)=xy=T_P(x,y)$$ and $$F(x,y)=x^2y.$$
\item [(\romannumeral 4)] For some $\beta<0$,
\begin{equation*}\label{eq120}
T(x,y)=\begin{cases}
(x^{\beta}+y^{\beta}-1)^{\frac{1}{\beta}} & \emph{if }(x,y)\in(0,1]^2,\\0 & \emph{otherwise}
\end{cases}
\end{equation*}
and $$F(x,y)=\begin{cases}
  (x^{\beta}+(xy)^{\beta}-1)^{\frac{1}{\beta}} & \emph{if }(x,y)\in(0,1]^2,\\0 & \emph{otherwise}.
  \end{cases}$$
\item [(\romannumeral 5)] There exists $0<c<1$ such that
$$T(x,y)=\begin{cases}
0  & \emph{if }(x,y)\in(0,1)^2\setminus[c,1)^2,\\  \min\{x,y\}  & \emph{otherwise}
  \end{cases}$$
and
$$F(x,y)=\begin{cases}
0  & \emph{if }(x,xy)\in(0,1)^2\setminus[c,1)^2,\\  xy  & \emph{otherwise}.
  \end{cases}$$
\item [(\romannumeral 6)] $$T=T_D$$
and
$$F(x,y)=\begin{cases}
  0  & x\in[0,1),\\ y  & x=1.
\end{cases}
$$
\end{enumerate}
\end{theorem}

\section{Conclusions}

 In this article, we introduced the concept of a general pseudo-homogeneous t-norm and established the representation theorem of a general pseudo-homogeneous t-norm (see Definition \ref{def2.4} and Theorem \ref{th0003.1}, respectively). Comparing with the definition of a pseudo-homogeneous t-norm in \cite{Xie2016}, the function $F$ related to a general pseudo-homogeneous t-norm needn't be a continuous function that fulfills the boundary condition $F(x,1)=0\Leftrightarrow x=0$. So that we generalized the concept of a pseudo-homogeneous t-norm in \cite{Xie2016} greatly. In the future work, we shall consider a general pseudo-homogeneity of an aggregation function.


\begin{thebibliography}{9}
\bibitem{Alsina2006} C. Alsina, M. J. Frank, B. Schweizer, Associative Functions: Triangular Norms and Copulas, World Scientific Press, Singapore, 2006.
\bibitem{Cabrera2014} G. Cabrera, M. Ehrgott, A. Mason, A. Philpott, Multi-objective optimisation of positively homogeneous functions and an application in radiation therapy, Oper. Res. Lett. 42 (2014) 268-272.
\bibitem{Ebanks1998} B. Ebanks, Quasi-homogeneous associative functions, Int. J. Math. Math. Sci. 21 (1998) 351-358.
\bibitem{Fuss1978} M. Fuss, D. McFadden, Y. Mundlak, A Survey of Functional Forms in the Economic Analysis of Production, North-Holland, Amsterdam, 1978.
\bibitem{Jurio2013} A. Jurio, H. Bustince, M. Pagola, A. Pradera, R. R. Yager, Some properties of overlap and grouping functions and their application to image thresholding, Fuzzy Sets Syst. 229 (2013) 69-90.
\bibitem{Jurio2014} A. Jurio, H. Bustince, M. Pagola, P. Couto, W. Pedrycz, New measures of homogeneity for image processing: an application to fingerprint segmentation, Soft Comput. 18 (2014) 1055-1066.
\bibitem{Klement2000} E. P. Klement, R. Mesiar, E. Pap,  Triangular norms, Kluwer Academic Publisher, Dordrecht, 2000.
\bibitem{Kuczma2009} M. Kuczma, An Introduction to the Theory of Functional Equations and Inequalities: Cauchy's Equation and Jensen's Inequality, Birkh\"{a}user Verlag AG, Basel$\cdot$Boston$\cdot$Berlin, 2009.
\bibitem{T2005} T. R\"{u}ckschlossov\'{a}, Homogeneous aggregation operators, computational intelligence, theory and applications, bernd reusch, Adv. Soft Comput. 33 (2005) 555-563.
\bibitem{Xie2016} A. F. Xie, Y. Su, H. W. Liu, On pseudo-homogeneous triangular norms, triangular conorms and proper uninorms, Fuzzy Sets Syst. 287 (2016) 203-212.

\end{thebibliography}
\end{document}